# ON THE HARDNESS OF ALMOST ALL SUBSET SUM PROBLEMS BY ORDINARY BRANCH-AND-BOUND

**Asst. Prof. Dr. Mustafa Kemal Tural**

Department of Industrial Engineering, Middle East Technical University, tural@metu.edu.tr

**SUMMARY**
Given $n$ positive integers $a_1, a_2, \ldots, a_n$, and a positive integer right hand side $\beta$, we consider the feasibility version of the subset sum problem which is the problem of determining whether a subset of $\{a_1, a_2, \ldots, a_n\}$ adds up to $\beta$. We show that if the right-hand side $\beta$ is chosen as $\lfloor r \sum_{j=1}^{n} a_j \rfloor$ for a constant $0 < r < 1$ and if the $a_j$'s are independent and identically distributed from a discrete uniform distribution taking values $\{1, 2, \ldots, \lfloor 10^{n/2} \rfloor\}$, then the probability that the instance of the subset sum problem generated requires the creation of an exponential number of branch-and-bound nodes when one branches on the individual variables in any order goes to 1 as $n$ goes to infinity.

## INTRODUCTION

Given $n$ positive integers $a_1, a_2, \ldots, a_n$, and a positive integer right hand side $\beta$, the subset sum problem (SSP) is the problem of finding a subset of $\{a_1, a_2, \ldots, a_n\}$ whose elements add up to $\beta$. Introducing a binary variable $x_j$ for each $j \in \{1, 2, \ldots, n\}$ which takes the value 1 if and only if $a_j$ is included in the subset, the SSP can be formulated as the following integer programming feasibility problem

$$\begin{aligned} ax &= \beta \\ x &\in \{0,1\}^n. \end{aligned} \quad \text{(SUB)}$$

The SSP is an NP-complete problem (Garey and Johnson, 1979). In the literature, several public-key cryptosystems are proposed based on the SSP (Merkle and Hellman, 1978; Shamir, 1983). The density $d$ of an SSP instance is defined as $d = n / \log_2 \max\{a_1, a_2, \ldots, a_n\}$ (Lagarias and Odlyzko, 1985). In cryptographic applications, SSP instances of low density are of particular importance because (SUB) may have several solutions when the density is high.

Brickell (1984) and Lagarias and Odlyzko (1985) focus on solving feasible low density SSPs and they show that almost all feasible low density subset sum problems can be solved in polynomial time. We refer the reader to Coster et al., (1992) for an improvement of the method of Lagarias and Odlyzko (1985). On the other hand, Furst and Kannan (1989) look into both feasible and infeasible SSP instances and show that almost all (feasible and infeasible) low density SSPs can be solved in polynomial time. In a related study, Pataki, Tural, and Wong (2010) generalize the result of Furst and Kannan (1989) from subset sum problems to bounded integer programming problems. In all of these methods, lattice-based techniques are used.

The most commonly used exact algorithms to solve integer programming problems are branch-and-bound, cutting plane methods, and branch-and-cut which is a combination of branch-and-bound and cutting plane methods. In this paper, we show that almost all low density subset sum problems are hard for ordinary branch-and-bound. In particular, we show that if the right-hand side $\beta$ is chosen as $\lfloor r \sum_{j=1}^{n} a_j \rfloor$ for a constant $0 < r < 1$ and if the $a_j$'s are independent and identically distributed (iid) from a discrete uniform distribution taking values $\{1,2, \ldots, \lfloor 10^{n/2} \rfloor\}$ in (SUB), then the probability that the instance of the SSP generated requires the creation of an exponential number of branch-and-bound nodes when one branches on the individual variables in any order goes to 1 as $n$ goes to infinity.

Our result is built on a result of Chvátal (1980) who identified a class of instances of the $0 - 1$ knapsack problem that are difficult to solve by a class of algorithms, called recursive, that use branch-and-bound, dynamic programming, and rudimentary divisibility arguments. The problem considered by Chvátal (1980) is the following optimization problem

$$\begin{aligned} maximize \quad & ax \\ subject\ to \quad & ax \leq \left\lfloor \frac{\sum_{j=1}^{n} a_j}{2} \right\rfloor \\ & x \in \{0,1\}^n. \end{aligned} \quad (KP)$$

Note that (KP) is always feasible, but (SUB) can be feasible or infeasible. Our result is based on the observation that an overwhelming majority of the instances of (SUB) is infeasible and an overwhelming majority of the infeasible instances are hard to solve by ordinary branch-and-bound.

In a related study, Krishnamoorthy (2008) considers infeasible equality-constrained knapsack problems and derives lower and upper bounds on the number of branch-and-bound nodes enumerated in a branch-and-bound method.

**METHOD AND RESULTS**
In this section, we first state our main theorem and then prove it using some lemmas.

**Theorem 1.** Let $r$ and $\epsilon$ be real numbers between 0 and 1 and $\beta = \lfloor r \sum_{j=1}^{n} a_j \rfloor$. Let each $a_j$ be iid from a discrete uniform distribution taking values $\{1,2, \ldots, M\}$ in (SUB), where $M = \lfloor 10^{n/2} \rfloor$. Then the probability that the instance of the SSP generated requires the creation of at least $2^{n^{1-\epsilon}}$ branch-and-bound nodes when one branches on the individual variables in any order goes to 1 as $n$ goes to infinity.

The idea of the proof of Theorem 1 is taken from Chvátal (1980). We first fix a constant $k$ that satisfies $0 < k < \epsilon < 1$. We then show that the probability that the coefficients $a_j$ satisfy the following two properties goes to 1 as n goes to infinity:

**Property 1.** $\sum_{i \in I} a_i \leq \frac{1}{n^k} \sum_{j=1}^{n} a_j$ whenever $|I| \leq n^{1-\epsilon}$.
**Property 2.** There is no subset $I$ of $\{1,2, \ldots, n\}$ with $\sum_{i \in I} a_i = \lfloor r \sum_{j=1}^{n} a_j \rfloor$.

If coefficients $a_j$ satisfy Property 2, then it is clear that the corresponding SSP is infeasible.

**Lemma 1.** The probability that the coefficients $a_j$ satisfy Property 1 goes to 1 as $n$ goes to infinity.

**Proof of Lemma 1.** If Property 1 is violated, then there exists a subset $I$ of $\{1, 2, \ldots, n\}$ with $\sum_{i \in I} a_i > \frac{1}{n^k} \sum_{j=1}^{n} a_j$ and $|I| \leq n^{1-\epsilon}$. As each $a_j$ is less than or equal to $M$, we get that

$$\sum_{j=1}^{n} a_j < n^k (M n^{1-\epsilon}) = M n^{1+k-\epsilon}. \tag{1}$$

To get an explicit upper bound on the probability that Property 1 is violated, we use the following identity

$$\sum_{\substack{i \text{ integer:} \\ i \geq (p+t)n}} \binom{n}{i} p^i (1-p)^{n-i} < e^{-2t^2 n}, \tag{2}$$

which is valid for $0 < p < 1$ and $t \geq 0$. By (1), we have that at least $n - 2n^{1+k-\epsilon}$ of the $a_i$'s must be less than or equal to $M/2$, as otherwise $\sum_{j=1}^{n} a_j \geq (2n^{1+k-\epsilon}) M/2 = M n^{1+k-\epsilon}$. Taking $p = \lfloor M/2 \rfloor / M$ and $t = 1/2 - 2n^{k-\epsilon}$ in (2), we get that

$$\sum_{\substack{i \text{ integer:} \\ i \geq (\frac{1}{2}+t)n}} \binom{n}{i} p^i (1-p)^{n-i} = \sum_{\substack{i \text{ integer:} \\ i \geq n - 2n^{1+k-\epsilon}}} \binom{n}{i} p^i (1-p)^{n-i}$$

$$\leq \sum_{\substack{i \text{ integer:} \\ i \geq (p+t)n}} \binom{n}{i} p^i (1-p)^{n-i} < e^{-2t^2 n} = e^{-2n\left(\frac{1}{2} - 2n^{k-\epsilon}\right)^2}$$

which goes to 0 as $n$ goes to infinity. In the above expressions, the second sum is the probability that at least $n - 2n^{1+k-\epsilon}$ of the $a_i$'s are less than or equal to $M/2$. As this probability goes to 0 as $n$ goes to infinity, we get that the probability that the coefficients $a_j$ satisfy Property 1 goes to 1 as $n$ goes to infinity.

**Lemma 2.** The probability that the coefficients $a_j$ satisfy Property 2 goes to 1 as $n$ goes to infinity.

**Proof of Lemma 2.** Let us fix a $0-1$ vector $x_1^*, x_2^*, \ldots, x_n^*$ and find an upper bound on the number of $a$ vectors satisfying $\sum_{j=1}^{n} a_j x_j^* = \lfloor r \sum_{j=1}^{n} a_j \rfloor$. Each $a$ vector satisfying $\sum_{j=1}^{n} a_j x_j^* = \lfloor r \sum_{j=1}^{n} a_j \rfloor$ also satisfies $r \sum_{j=1}^{n} a_j - 1 < \sum_{j=1}^{n} a_j x_j^* \leq r \sum_{j=1}^{n} a_j$. Hence, we have that $r \sum_{j=1}^{n-1} a_j - 1 - \sum_{j=1}^{n-1} a_j x_j^* < a_n(x_n^* - r) \leq r \sum_{j=1}^{n-1} a_j - \sum_{j=1}^{n-1} a_j x_j^*$. This implies that once $a_1, a_2, \ldots, a_{n-1}$ are fixed, we have at most $k(r) = \lceil \max\{1/r, 1/(1-r)\} \rceil$ many choices for $a_n$. This follows from the two values $x_n^*$ may take. If $x_n^* = 0$, we have less than or equal to $\lceil 1/r \rceil$ many choices for $a_n$; and if $x_n^* = 1$, we have less than or equal to $\lceil 1/(1-r) \rceil$ many choices for $a_n$. Therefore for a fixed a $0-1$ vector $x_1^*, x_2^*, \ldots, x_n^*$, the number of $a$ vectors satisfying $\sum_{j=1}^{n} a_j x_j^* = \lfloor r \sum_{j=1}^{n} a_j \rfloor$ is at most $k(r) M^{n-1}$. As there are $2^n$ choices for $x_1^*, x_2^*, \ldots, x_n^*$, the probability that the coefficients $a_j$ violate Property 2 is at most $k(r) M^{n-1} 2^n / M^n = k(r) 2^n / M = k(r) 2^n / \lfloor 10^{n/2} \rfloor$ which goes to zero as $n$ goes to infinity.

**Lemma 3.** For positive coefficients $a_j$ satisfying Properties 1 and 2, if $\beta \in \left[\frac{1}{n^k} \sum_{j=1}^{n} a_j, \left(1 - \frac{1}{n^k}\right) \sum_{j=1}^{n} a_j\right]$ and if (SUB) is infeasible, then the ordinary branch-and-bound creates at least $2^{n^{1-\epsilon}}$ branch-and-bound nodes when one branches on the individual variables in any order.

**Proof of Lemma 3.** We will show that none of the nodes in the branch-and-bound tree is pruned by infeasibility unless more than $n^{1-\epsilon}$ of the variables are fixed.

Assume that $\leq n^{1-\epsilon}$ are fixed to 0 or 1. Let $I$ and $\bar{I}$ be the set of indices of the fixed and unfixed variables, respectively. By Property 1, we have that $\sum_{i \in I} a_i \leq \frac{1}{n^k} \sum_{j=1}^{n} a_j$ and $\sum_{i \in \bar{I}} a_i \geq \left(1 - \frac{1}{n^k}\right) \sum_{j=1}^{n} a_j$. Therefore, we get that $\sum_{i \in I} a_i \leq \frac{1}{n^k} \sum_{j=1}^{n} a_j \leq \lfloor r \sum_{j=1}^{n} a_j \rfloor = \beta \leq \left(1 - \frac{1}{n^k}\right) \sum_{j=1}^{n} a_j \leq \sum_{i \in \bar{I}} a_i$. Therefore, by assigning fractional values to $x_i, i \in \bar{I}$, we can obtain a feasible solution to the LP relaxation of (SUB). This is because if we assign the value 0 to all unfixed variables, a value, say $\alpha$, is obtained which is less than or equal to $\beta$. As assigning values in [0,1] to the unfixed variables, any value between $[0, \beta]$ can be obtained (because $\beta \leq \sum_{i \in \bar{I}} a_i$), we can obtain $\beta - \alpha$ as well.

**Proof of Theorem 1.** We are now ready to prove Theorem 1. Lemmas 1 and 2 imply that the probability that the coefficients $a_j$ satisfy Properties 1 and 2 goes to 1 as n goes to infinity. As $n$ goes to infinity, the probability that $\frac{1}{n^k} \leq r \leq 1 - \frac{1}{n^k}$ is satisfied goes to one. Therefore as n goes to infinity, the probability that $\beta \in \left[\frac{1}{n^k} \sum_{j=1}^{n} a_j, \left(1 - \frac{1}{n^k}\right) \sum_{j=1}^{n} a_j\right]$ goes to one. In other words, conditions in Lemma 3 are satisfied with a probability that goes to one as $n$ goes to infinity. Finally, whenever Lemma 3 is satisfied, at least $2^{n^{1-\epsilon}}$ branch-and-bound nodes are created. So, we have proved that the probability that the instance of the SSP generated requires the creation of at least $2^{n^{1-\epsilon}}$ branch-and-bound nodes when one branches on the individual variables in any order goes to 1 as $n$ goes to infinity.

## CONCLUSION
In this paper, we have shown that an overwhelming majority of the low density subset sum problems are hard to solve by ordinary branch-and-bound. This result complements the positive results on the solvability of the majority of the low density subset problems. In our proof, we have first argued that almost all subset sum problems generated in accordance with our generation procedure are infeasible. We have then proved that for almost all such infeasible subset sum problems an exponential number of branch-and-bound nodes are created.